\documentclass[11pt, a4paper]{article}  
\author{Frank Schuhmacher}

\title{Three deformation functors for associative algebras}

\makeindex 

\usepackage{mathenv}
\usepackage[intlimits]{amsmath}
\usepackage[T1]{fontenc}
\usepackage{amssymb}
\usepackage{amsthm} 
\usepackage{amscd}
\usepackage[all]{xypic}
\usepackage{makeidx}

\newcounter{punkt}

\theoremstyle{definition}
\newtheorem{defi}{Definition}[section]

\newtheorem{beisp}[defi]{Example}

\theoremstyle{theorem}
\newtheorem{satz}[defi]{Theorem}

\newtheorem{lemma}[defi]{Lemma}
\newtheorem{kor}[defi]{Corollary}
\newtheorem{prop}[defi]{Proposition}

\newcommand{\nach}{\longrightarrow}
\newcommand{\sub}{\subseteq}

\newcommand{\isom}{\cong}

\newcommand{\ZZ}{\mathbb{Z}}

\newcommand{\ppp}{\cdot\ldots\cdot}
\newcommand{\kmod}{k\text{-Mod}}
\newcommand{\Bmod}{B\text{-Mod}}

\newcommand{\Balg}{B\text{-Alg}}

\newcommand{\B}{\mathcal{B}}

\newcommand{\sX}{\mathcal{X}}

\renewcommand{\Im}{\operatorname{Im}}

\newcommand{\m}{\mathfrak{m}}

\newcommand{\spec}{\operatorname{Spec}}

\newcommand{\Hom}{\operatorname{Hom}}

\newcommand{\Tor}{\operatorname{Tor}}

\newcommand{\Kern}{\operatorname{Kern} }

\newcommand{\ot}{\otimes}

\newcommand{\lie}{[\cdot,\cdot]}

\newcommand{\Coder}{\operatorname{Coder}}

\newcommand{\ots}{\otimes\ldots\otimes}

\newcommand{\pnu}{\downarrow}

\newcommand{\pno}{\uparrow}

\newcommand{\Def}{\operatorname{Def}}

\newcommand{\Der}{\operatorname{Der}}

\newcommand{\obj}{\operatorname{Ob}}

\newcommand{\DGCoalg}{\operatorname{DG-Coalg}}
\newcommand{\DGAlg}{\operatorname{DG-Alg}}

\newcommand{\Aut}{\operatorname{Aut}}
\newcommand{\End}{\operatorname{End}}
\newcommand{\ass}{\operatorname{ass}}
\newcommand{\rel}{\operatorname{rel}}
\renewcommand{\flat}{\operatorname{flat}}
\newcommand{\MC}{\mathcal{MC}}

\newcommand{\Art}{\mathfrak{Art}}
\newcommand{\set}{\mathfrak{set}}
\newcommand{\ad}{\operatorname{ad}}

\begin{document}
\maketitle

\begin{abstract}
We show that three deformation functors
(deformations of the product, flat deformations and
deformations of the relations)
assigned to an associative algebra are naturally
isomorphic. 
\end{abstract}

\section*{Introduction}

In the litterature (see \cite{Gerst1} and references here-in),
deformations of an associative algebra $A$ mostly refer to
deformations of the associative multiplication of $A$.
But, dually to deformations of schemes and singularities, it
is as natural to consider flat deformations of
associative algebras. A deformation of the multiplication
of $A$ defines directly a flat deformation. In this paper
we show that, up to equivalence, any flat deformation
over an Artin base is just a deformation of the 
multiplication.

A third approach is to consider $A$ as a quotient of a free
algebra, modulo an ideal of relations, and to study deformations
of those relations. To make this approach fit into the flat
concept, we also have to control ``higher relations'' which leads
to deforming the differential of a whole free DG algebra resolution
of $A$. By abuse of language, we still speak of
deformations of the relations of $A$. 
It was shown by Kontsevich \cite{KontL} that deformations of
the associative multiplications of $A$, parametrized by an
Artin algebra $B$ (or by a projective limit of Artin algebras), 
correspond one-to-one to deformations of the relations of
$A$, parametrized by $B$. 
The aim of this paper is
to show that flat deformations over Artin algebras 
are just a third perspective to ``deformations of $A$''.\\

\noindent
\textbf{Conventions:}
For the whole paper, fix a ground field
$k$ of characteristic zero. We will study deformations
of an associative non-unitary $k$-algebra $A$.
All constructions and statements are true with sleight modifications
for unitary algebras. 
The multiplication $A\ot_k A\nach A$ will be denoted by $\alpha$.
Deformations will be parametrized by a commutative
local artinian algebra
$B=(B,\m)$ such that $B/\m=k$. Let $\Art$ be the category of these algebras.
For any $k$-module (or DG module) $M$, set
$M_B:=M\ot_kB$ and $M_\m:=M\ot_k\m$.

\section{Deformation theory via DG Lie algebras}\label{secti}

A good introduction to deformation theory via DG Lie algebras is
given in \cite{Mane}. We will just recall the general basic
concept: To study deformations of a given object $A$ in some category,
we assign a deformation functor $\Def$ to $A$,
i.e. the set-valued
functor $\Def:\Art\nach\set$ on the category of Artin algebras,
assigning to $B\in\obj(\Art)$ the quotient of the
set of all deformations of $A$ parametrized by $B$, modulo
equivalent deformations. We say that the deformations of $A$ are 
\textbf{governed by} a DG Lie algebra $L=(L,d,\lie)$, 
if the functor $\Def$ is of the form
\begin{equation}\label{defmc}
\Def(B)=\MC(B)/\sim,
\end{equation}
where the right hand-side is defined as follows:
The tensor products $L_B:=L\ot B$ and $L_\m:=L\ot\m$
are again DG Lie algebras with structure maps
defined by $d(u\ot b)=du\ot b$ and
$[u\ot b,u'\ot b']=[u,u']\ot bb'$. Since the ideal $\m$ is
nilpotent, the Lie algebra $L^0_\m$ is also nilpotent.
We can thus form the group $G:=\exp(L^0_\m)$, and the exponential
map $\exp:L^0_\m\nach G$ is bijective.
The group $G$ operates on $L^1_\m$ via
$$\exp(f)\cdot u=\exp(\ad_f)(u),$$
where $\ad_f$ is the endomorphism $u\mapsto [f,u]$ of $L^1_\m$.
It is convenient to introduce the DGL algebra
$\tilde{L}_\m=(\tilde{L}_\m,d_d,\lie_d)$, 
defined by $\tilde{L}_\m^i:=L_\m^i$, for
$i\neq 1$ and  $\tilde{L}_\m^1:=L_\m^1\oplus k\cdot d$.
The structure maps of $\tilde{L}_\m$ are defined
by $d_d(u+\lambda d)=du$, $[d,u]_d=du$ and $[d,d]_d=0$, for
$\lambda\in k$ and $u\in L_\m$.
We have the affine map $L_\m\nach\tilde{L}_\m$, given by
$u\mapsto\tilde{u}:=u+d$.
Since $\tilde{L}_\m^0=L^0_\m$, the group $G$ also acts on
$\tilde{L}^1_\m$. Now we can define $\MC(B)$ which is the
set
$$\{u\in L^1_\m|\;du+\frac12[u,u]=0\}$$
of \textbf{Maurer-Cartan-elements} of $L_\m$.
The tilde-map restrics to a bijection between
$\MC(B)$ and the set
$$\tilde{MC}(B):=
\{\tilde{u}\in\tilde{L}^1_\m|\;[\tilde{u},\tilde{u}]_d=0\}.$$
The set $\tilde{MC}(B)$ is invariant under the action of $G$.
We call two Maurer-Cartan elements $u,v\in\MC(B)$
\textbf{(gauge-)equivalent}, if $\tilde{u}$ and
$\tilde{v}$ are in the same orbit of the action of
$G$ on $\tilde{MC}(B)$.
This defines the equivalence relation in equation~\ref{defmc}.\\

We can verify immediately that if two deformation
theories are governed by DG Lie algebras $L_a$ and $L_b$,
respectively, and if there is an (iso)morphism
$L_a\nach L_b$ of DG Lie algebras, we get an
(iso)morphism of deformation functors. But due to Kontsevich \cite{Kont},
there is a much stronger statement:

\begin{satz}\label{isodef}
A quasi-isomorphism between two DG Lie algebras 
(in fact a quasi-isomorphism in the category of $L_\infty$-algebras
is sufficient) induces
an isomorphism between the corresponding deformation
functors. 
\end{satz}

Below, to any associative algebra,
we will assign two deformation functors
$\Def_{\rel}$ (deformations of relations)
and $\Def_{\ass}$ (associative deformations)
that are governed by DG Lie algebras 
$L_{\ass}$ and $L_{\rel}$, respectively,
and a deformation functor $\Def_{\flat}$, for which it is not
clear, a priori, if it is governed by
a DG Lie algebra.
It was shown by Kontsevich
\cite{KontL} that the DG Lie algebras $L_{\rel}$ and $L_{\ass}$ are
quasi-isomorphic, thus the deformation functors
$\Def_{\ass}$ and $\Def_{\rel}$ are isomorphic. The aim of this
paper is to show that the third one $\Def_{\flat}$ is
isomorphic to the both of them.

\section{Deformations of the associative structure}\label{Sass}

\begin{defi}\label{defass}
A \textbf{deformation of the associative structure}
of $A$ with basis $B=(B,\m)$ is a $B$-linear map
$\beta:A_B\ot_BA_B\nach A_B$ with $\Im(\beta)\sub A_\m$
such that 
the sum $\alpha+\beta$ is associative
as a product on $A_B$.
\end{defi}

\begin{defi}\label{equivalent}
Two deformations $\beta$ and $\beta'$ of the associative
structure of $A$ with basis $B$
are \textbf{equivalent}, if there exists a
$B$-linear isomorphism 
$\Phi:(A_B,\alpha+\beta)\nach (A_B,\alpha+\beta')$ 
(in the unitary case, $\Phi$ has to be unitary)
such that
the following diagram commutes:
$$\xymatrix{
 & A &&\\
A_B\ar[ur]\ar[rr]^{\Phi}&&A_B\ar[ul]& (\ast)}$$
Maps without label are canonical.
\end{defi}

It was discovered by Gerstenhaber, that the
deformations of the associative structure of $A$ are
governed by a DG Lie algebra $L_{\ass}$, namely the Hochschild
cochain complex. The graded Lie bracket on this complex is
now called the ``Gerstenhaber bracket''. Let's recall
the modern definition of the Hochschild cochain complex:

In the sequel, $M$ shall always denote
the graded $k$-module $\pnu A$ which is defined
as $(\pnu A)^{-1}=A$ and $(\pnu A)^i=0$, for $i\neq -1$.
We have natural maps $\pnu:A\nach M$ (degree -1) and
$\pno:M\nach A$ (degree +1).
For any graded $k$-module $W$, let
$TW:=\coprod_{n\geq 1}W^{\ot n}$ be the (non-unitary)
tensor algebra. On $TW$, we have a coproduct
$\Delta:TW\nach TW\ot TW$, defined by
$$w_1\ots w_n\mapsto\sum_{i=1}^{n-1}(w_1\ots w_i)\ot(w_{i+1}\ots w_n).$$
Recall that the \textbf{bar construction} $\B A$ of $A$ is the 
DG coalgebra $(TM,\Delta,Q)$, where the degree one codifferential
$Q$ on $TM$ is uniquely defined by its (co)restriction
\begin{align*}
Q_2:M\ot M&\nach M\\
m_1\ot m_2&\mapsto -\pnu\alpha(\pno m_1,\pno m_2).
\end{align*}
For more details, see \cite{Smir}.
The \textbf{Hochschild cochain complex} of $A$
is the DG Lie algebra $L_{\ass}:=\Coder(\B A)=\coprod_i\Hom(M^{\ot i},M)$
of coderivations on the DG coalgebra $\B A$. Its Lie bracket
is the graded commutator and the (degree one) differential $d_{\ass}$ is
given by $q\mapsto [Q,q]$.
Observe that
$$L_{\ass}^0(\m)\isom\Hom_{\kmod}(A,A_\m)\isom
\{f\in\End_{\Bmod}(A_B)|\;\Im f\sub A_\m\}.$$
Now, $L^0_{\ass}(\m)$
is a nilpotent 
Lie subalgebra of $\End_{\Bmod}(A_B)$. The Lie group $G_{\ass}$ 
associated to $L_{\ass}(\m)^0$ is the subgroup
$$\{\phi\in\Aut_{\Bmod}(A_B)|\;\phi=1+\psi\text{ for a }\psi\text{ with }
\Im \psi\sub A_{\m} \}$$ of $\Aut_{\Bmod}(A_B)$.
The exponetial map
\begin{align*}
\exp:L_{\ass}(\m)^0&\nach G_{\ass}\\
f&\mapsto 1+ f+\frac 12 f^2+\ldots
\end{align*}
is bijective.
Its inverse is given by $1+\psi\mapsto\ln(1+\psi)=\sum_{n=1}^\infty
\frac{(-1)^{n-1}\psi^n}{n}$. 
Observe that
\begin{align*}
L^1_{\ass}(\m)&\isom\Hom_{\kmod}(A\ot_kA,A_\m)\\
&\isom
\{\beta\in\Hom_{\Bmod}(A_B\ot_BA_B,A_B)|\;\Im\beta\sub A_\m\}.
\end{align*}

The group $G_{\ass}$ acts on $L_{\ass}^1(B)$
by $\phi:\gamma\mapsto\phi\circ\gamma\circ(\phi^{-1}\ot\phi^{-1})$.
This is just the group action described more abstractly in
Section~\ref{secti}.
We identify deformations $\beta$ and $\beta'$ of $\alpha$,
if there exists an element $\phi$ of $G_{\ass}$ such that
$\phi.(\alpha+\beta)=\alpha+\beta'$. (Be careful, in this case, 
we don't have $\beta'=\phi.\beta$!)
This is just the equivalence relation in Definition~\ref{defass}.
Verify that
deformations of the
associative structure of $A$ with basis $B$ correspond
to Maurer-Cartan elements 
$$\MC_{\ass}(B):=\{\beta\in L^1_{\ass}(\m)|\;
d_{\ass}(\beta)+\frac{1}{2}[\beta,\beta]=0\}.$$
Each element $\phi$ of $G_{\ass}$, defines a bijection  
(that we shall also denote by $\phi$)
\begin{align*}
\phi:\MC_{\ass}(B)&\nach\MC_{\ass}(B)\\
\beta&\mapsto\phi.(\alpha+\beta)-\alpha.
\end{align*}
Two elements $\beta$ and $\beta'$ in $\MC_{\ass}(B)$ 
are equivalent, i.e. if $\beta'=\phi(\beta)$, for a $\phi\in G_{\ass}$, 
if and only if they correspond to equivalent deformations
of $\alpha$. 
We define the deformation functor 
$\Def_{\ass}:\Art\nach\set$ as
$$\Def_{\ass}(B):=\MC_{\ass}(B)/\sim.$$

\section{Deformations of the relations}

We can choose a free non-positively graded
DGA resolution $(R,s)$ of $A$, i.e. $(R,s)$ is 
the free graded tensor algebra over a graded $k$-module
together with a derivation $s$ of degree $+1$
with $s^2=0$ such that $(R,s)$ is acyclic in negative 
degrees and $H^0(R,s)\isom A$ via a fixed morphism
$R^0\nach A$. It can easily be verified that such a resolution
exists and we will see later that there is even a
canonical one.

\begin{defi}
A \textbf{deformation of the relations} of $A$ with basis $B$ 
is a degree one derivation $\delta$ of the free graded $B$-algebra
$R_B=R\ot_kB$ such that the following two conditions hold:
\begin{enumerate}
\item
$\Im\delta\in R\ot_k\m$,
\item
$(s+\delta)^2=0$, i.e. the perturbed derivation $s+\delta$
on $R_B$ is again a differential.
\end{enumerate} 
\end{defi} 
   
First, we want to
prove that, for such a deformation $\delta$, the complex
$(R_B,s+\delta)$ is acyclic in negative degrees.

\begin{prop}\label{specs}
For each deformation $\delta$ of the relations of
$A$ with base $B$, the complex $(R_B,s+\delta)$ is a free
resolution of $\tilde A:=H^0(R_B,s+\delta)$ over $B$. 
\end{prop}
\begin{proof}
The inclusions $...\sub\m^2\sub\m\sub B$ give raise to a filtration
of the complex $(R_B,s+\delta)$.
The term $E_1^{p,q}$ of the spectral sequence assigned to
this filtration  is zero, 
for $p+q\neq 0$. The statement is thus a consequence of the
Classical Convergence Theorem (see \cite{Weib}) for
spectral sequences. 
\end{proof}

\begin{defi}
Two Deformations $\delta_1$ and $\delta_2$ of the relations of $A$
with basis $B$ are \textbf{equivalent}, if there exist 
a B-linear isomorphism 
$\Phi:(R_B,s+\delta_1)\nach (R_B,s+\delta_2)$ 
of DG algebras such that the diagram
$$\xymatrix{
&R&&\\
R_B\ar[ur]\ar[rr]^{\Phi} && R_B\ar[ul] & (\ast\ast)
}$$
commutes.
\end{defi}

Again, we will construct a DG Lie algebra, governing the
deformations of the relations of $A$.
Let $L_{\rel}$ be the graded module $\Der(R)=\coprod_{i\in\ZZ}\Der^i(R)$
of $k$-linear derivations of $R$. It carries a DGL structure,
where the Lie bracket is again the graded commutator and
the degree one differential $d_{\rel}$
is given by $\delta\mapsto [s,\delta]$.
The complex $(L_{\rel},d_{\rel})$ is called \textbf{
(noncommutative) tangent complex}.
As above, consider the DG Lie algebras $L_{\rel}(B):=L_{\rel}\ot B$ and
$L_{\rel}(\m):=L_{\rel}\ot \m$. The Liegroup corresponding to
the nilpotent Lie algebra $L^0_{\rel}(\m)$ is given by
$$G_{\rel}=\{\phi\in\Aut_{\Balg}(A_B)|\;\phi=1+\psi\text{ for a }
\psi\text{ with }\Im\psi\sub A_\m\}.$$
Remark that in contrast to $G_{\ass}$, the elements of $G_{\rel}$
are algebra homomorphisms.    The exponential map
\begin{align*}
\exp:L_{\rel}(\m)^0&\nach G_{\rel}\\
f&\mapsto 1+f+\frac 12 f^2+\ldots
\end{align*}
is bijective.
The group $G_{\rel}$ acts on $L^1_{\rel}(B)$ by
$\phi:\gamma\mapsto\phi\gamma\phi^{-1}$.
We identify deformations $\delta$ and $\delta'$ of the relations of $A$ with
basis $B$, if there exists an element of $G_{\rel}$ such that
$\phi.(s+\delta)=s+\delta'$.
Verify that the
deformations of the relations of $A$ are in one-to-one
correspondence with the Maurer-Cartan elements:
$$\MC_{\rel}(B):=\{\delta\in L_{\rel}(\m)^1|\;
d_{\rel}(\delta)+\frac{1}{2}[\delta,\delta]=0 \}$$ 
Each element $\phi$ of $G_{\rel}$, defines a bijection  
(that we shall also denote by $\phi$)
\begin{align*}
\phi:\MC_{\ass}(B)&\nach\MC_{\ass}(B)\\
\beta&\mapsto\phi.(\alpha+\beta)-\alpha.
\end{align*}
Two elements $\delta$ and $\delta'$ in $\MC_{\rel}(B)$ 
are equivalent, i.e. if $\delta'=\phi(\delta)$, for a $\phi\in G_{\rel}$, 
if and only if they correspond to equivalent deformations
of the relations of $A$ with basis $B$. 
We define the deformation functor 
$\Def_{\rel}:\Art\nach\set$ as
$$\Def_{\rel}(B):=\MC_{\rel}(B)/\sim,$$
where $\sim$ denotes the relation generated by this group action.
The following theorem, taken from \cite{KontL}, proves that
the deformation functor does not depend on the
chosen resolution $R$ of $A$:

\begin{satz}
The tangent Lie algebras $L_{\rel,1}$ and $L_{\rel,2}$,
corresponding to  
two different free DGA resolutions $R_1$ and $R_2$ 
of $A$ are quasi-isomorphic as DGL algebras. As consequence,
they define isomorphic deformation functors.
\end{satz}

\section{Flat deformations}

\begin{defi}
A \textbf{flat deformation} of $A$ with basis $B$ is
a flat associative $B$-algebra $\tilde A$, 
together with a fixed isomorphism
$\tilde A\ot_Bk\nach A$.
\end{defi} 

\begin{beisp}
For $A$ commutative, the
deformations of the affine scheme $X=\spec A$
over $\spec B$ are by definition cartesian diagrams
$$\xymatrix{
\sX\ar[d]^g & X\ar[d]\ar[l]\\
\spec B & \ast\ar[l]}$$
If  $\sX=\spec\tilde A$ is an affine flat scheme over $\spec B$,
then $\tilde A$ is a flat deformation of $A$.
\end{beisp}

\begin{defi}
Two flat deformations $\tilde A_1$ and $\tilde A_2$
of $A$ with basis $B$ are \textbf{equivalent}, if there
exists a $B$-algebra isomorphism
$\Phi:\tilde A_1\nach\tilde A_2$
such that the diagram
$$\xymatrix{
 & A &&\\
\tilde A_1\ar[ur]\ar[rr]^{\Phi} & & \tilde A_2\ar[ul]& (\ast\ast\ast)
}$$ 
commutes.
\end{defi}

Flat deformations of $A$ are not directly related to a
DG Lie algebra. We define the deformation functor
$\Def_{\flat}:\Art\nach\set$ by
$$\Def_{\flat}(B):=\{\text{flat deformations of $A$ with basis $B$}\}/\sim.$$

\subsection{From associative deformations to flat deformations}

If $\beta$ is a deformation of the associative structure
of $A$ with basis $B$, then $\tilde A:=(A_B,\alpha+\beta)$
is a flat deformation of $A$. If $\beta$ is equivalent to $\beta'$
as deformation of $\alpha$,
then  $(A_B,\alpha+\beta)$ is equivalent to $(A_B,\alpha+\beta')$ 
as flat deformation.
Thus we have a natural transformation
$$F:\Def_{\ass}\nach\Def_{\flat}$$
of deformation functors. By construction, for any
artinian basis $B$, the map
$F(B):\Def_{\ass}(B)\nach\Def_{\flat}(B)$
is injective.

\subsection{From deformations of the relations to flat deformations}

We fix a free DGA resolution $(R,s)$ of $A$. This data
contains an algebra map $\phi:R^0\nach A$, inducing an 
isomorphism $H^0(R,s)\nach A$.
Let $\delta$ be a deformations of the relations of $A$ with
basis $B$. We have seen that the DGA algebra $(R_B,s+\delta)$ is
a free resolution of $\tilde A:=H^0(R_B,s+\delta)$ over
$B$. 
\begin{prop}
The algebra $\tilde A$ is flat as $B$-module.
\end{prop}
\begin{proof}
For each $B$-module $M$ and each $i>0$, we have
$$\Tor_i^B(\tilde A,M)=H^{-i}(R_B\ot_BM)=0,$$
which follows by the same specatral sequence argument as in
the proof of Proposition~\ref{specs}.
\end{proof}
 
The assignment $\delta\mapsto H^0(R_B,s+\delta)$
defines a natural transformation
$$G:\Def_{\rel}\nach\Def_{\flat}$$ of deformation
functors. The naxt aim is to show that $G$ is surjective.
We will make use of Nakayama's Lemma for Artin algebras
(that works without any finitely generated assumption).
For the proof, see \cite{Mane}.

\begin{lemma}\label{nakayama}
For a local artinian algebra $(B,\m)$ and any $B$-module,
if $M\ot_BB/\m=0$, then $M=0$.
\end{lemma}

In the sequel, for any $B$-module $M$, we will use
the abreviation $M':=M\ot_BB/\m$.

\begin{lemma}\label{flaka}
Consider a complex
$$\xymatrix{...\ar[r] & M^{-2}\ar[r]^{d^{-2}} & M^{-1}\ar[r]^{d^{-1}} & 
M^0}$$
of flat $B$-modules such that the induced complex
$(M',d')$ is acyclic (in degrees $\leq i$). Then, $(M,d)$ is acyclic
(in degrees$ \leq i$ and $\Kern(d^{i-1})$ is flat).
The same statement is true, if we replace $B$ by any local noetherian 
algebra under the additional assumption that each $M^i$ is a
finitely generated $B$-module. 
\end{lemma}
\begin{proof}
First, we show that $d^{-1}$ is surjective:
We have a short exact sequence
$$0\nach \Im(d^{-1})\nach M^0\nach H^0(M)\nach 0.$$
Application of the right exact functor $\ot_BB/\m$ gives
an exact sequence
\begin{equation}\label{seq11}
\Im(d^{-1})'\nach {M'}^0\nach H^0(M)'\nach 0.
\end{equation}
Since the epimorphism
${d'}^{-1}$ factors through
$\Im(d^{-1})'$, the first arrow of the
exact sequence~\ref{seq11} is surjective, hence $H^0(M)'=0$.
By Lemma~\ref{nakayama}, we get $H^0(M)=0$.
Next we show that $H^{-1}(M)=0$:
Since we have an exact sequence
\begin{equation}\label{seq12}
0\nach \Kern(d^{-1})\nach M^{-1}\nach M^{0}\nach 0,
\end{equation}
the flatness of $M^0$ and $M^{-1}$ implies the flatness
of $\Kern(d^{-1})$. Applying th functor $\ot_BB/\m$ to the
exact sequence~\ref{seq12}, we get an exact sequence
\begin{equation}
0\nach \Kern(d^{-1})'\nach {M'}^{-1}\nach {M'}^{0}\nach 0,
\end{equation}
hence we have a natural isomorphism
$\Kern({d'}^{-1})\isom \Kern(d^{-1})'$.
Applying $\ot_BB/\m$ to the short exact sequence
$$0\nach \Im(d^{-2})\nach\Kern(d^{-1})\nach H^{-1}(M)\nach 0,$$
we get an exact sequence
$$\Im(d^{-2})'\nach\Kern(d^{-1})'\nach H^{-1}(M)'\nach 0.$$
Since the epimorphism
${M'}^{-2}\nach\Kern({d'}^{-1})$ factors through
$\Im(d^{-2})'$, we see that the first map of the last
exact sequence is surjective. Thus $H^{-1}(M)'=0.$
Again, by Lemma~\ref{nakayama}, $H^{-1}(M)=0$.
Continuing inductively in the same manner, we get that
for each $i\leq 1$, the module $\Kern(d^i)$ is flat over $B$ and that
$H^i(M)=0$.
\end{proof}

\begin{prop}
For each algebra $B$ in $\Art$, the map
$G(B):\Def_{\rel}(B)\nach\Def_{\flat}(B)$ is surjective.
\end{prop}
\begin{proof}
Let $f:B\nach\tilde A$ be a flat deformation of $A$ with basis
$B$. As graded algebra, set $\tilde R:=R_B$. 
Remark that $\tilde{R}'=\tilde R\ot_Bk=R$.\\

We claim that we can find (1)
a lifting $\tilde \phi:\tilde{R}^0\nach\tilde A$
of $\phi$, i.e. map of $B$-algebras such that
$\tilde{\phi}'=\phi$, and (2) a differential
$\tilde s$ on the graded algebra $\tilde R$ such that
$\tilde{s}'=s$ and such that
$\tilde\phi$ induces an isomorphism
$H^0(\tilde R,\tilde s)\nach\tilde A$.\\

Since $\tilde{R}^0$ is a free $B$-algebra, a lifting
$\tilde\phi$ of $\phi$ can be chosen. By Lemma~\ref{flaka}, this
lifting $\tilde\phi$ is surjective and thus, $\Kern\tilde\phi$
is flat over $B$. If follows that $(\Kern\tilde\phi)'\isom\Kern\phi$.
We define the differential $\tilde s$ inductively, by chosing
its values on the free generators $x$ of $R$.
For be a free generator $x$ of degree $-1$, there is a lifting
$\tilde r$ of $s(x)$ in $\Kern\tilde\phi$. Set $\tilde s(x):=\tilde r$.
In this way, we define the component 
$\tilde{s}^{-1}:\tilde{R}^{-1}\nach\tilde{R}^0$.
Applying Nakayama, as in the proof of Lemma~\ref{flaka}, we see that
$\Im(\tilde{s}^{-1})=\Kern\tilde\phi$, thus we have a short
exact sequence
$$0\nach\Kern\tilde{s}^{-1}\nach\tilde{R}^{-1}\nach\Kern\tilde\phi\nach 0.$$
It follows that $\Kern\tilde{s}^{-1}$ is flat over $B$ and that
$\Kern s^{-1}=(\Kern\tilde{s}^{-1})'$.
Now, for free generators $x$ of $R$ of degree $-2$,
we can define $\tilde{s}^{-2}(x)$ as a lift of $s^{-2}(x)$ in
$\Kern\tilde{s}^{-1}$ and so on.
Inductively, for each $n>0$, we can define the values
$\tilde{s}^{-n}(x)$, for free generators of $R$ of degree $-n$
in such a way that $\Kern\tilde{s}^{-n}$ is flat over $B$
and such that the complex
$$\Kern\tilde{s}^{-n}\nach \tilde{R}^{-n}\nach\cdots\nach\tilde{R}^0$$
is a resolution of $\tilde A$.
Thus the differential $\tilde s$ can be constructed
as desired.\\

Finally, we see that $\tilde s$ equals $s+\delta$, for a deformation
$\delta$ of the relations of $A$ with basis $B$.
\end{proof}

\section{The canonical free DGA resolution}

Consider a differential graded coalgebra
$C$ with comultiplication $\Delta$
and degree one codifferential $q$.
Set $W$ to be the shifted module $C[-1]=\pno C$.
Recall that the \textbf{cobar construction}
$\Omega C$ of $C$ is the free graded tensor
algebra 
$TW$ together with the degree one differential $d$,
which is uniquely defined by its corestrictions
$d_1:=\pno q\pnu:W\nach W$ and $d_2:=\pno^{\ot 2}\Delta\pnu$.
The cobar construction is a functor
$$\Omega:\DGAlg\nach\DGCoalg.$$
The following statement can be found in \cite{Lefe} or \cite{Smir}.

\begin{satz}\label{brat}
\begin{enumerate}
\item
Bar and cobar construction form an adjoint pair
$(\Omega,\B)$ of functors. (With appropriate
model structures on the categories $\DGAlg$ and $\DGCoalg$
they even form a Quillen-functor.)
\item
For any DG algebra $A$,
the natural map $\Omega \B A\nach A$ is a quasi-isomorphism.
I.e. $\Omega \B A$ is a canonical free DGA resolution of $A$.
\end{enumerate}
\end{satz}

It is worth, studying more closely the canonical
free DGA resolution $\Omega\B A$ of an associative
$k$-algebra $A$. As earlier, set $M:=\pnu A$.
Recall that the bar construction $\B A$ is the DG coalgebra $TM$ with
coalgebra structure as in Section~\ref{Sass}. 
As graded object,
the complex $\Omega\B A$ equals $T\pno(TM)$.
Picturally:

$$\xymatrix{
 & & & \pno(M\ot M\ot M\ot M)\ar[dl]\ar[d]\\
 & & \pno(M\ot M\ot M)\ar[dl]\ar[d] & 
\text{4 factors,}\atop\text{2 arrows}\ar[dl]\ar[d]\\
 & \pno(M\ot M)\ar[dl]\ar[d] & \pno(M\ot M)\ot\pno M\atop
+\pno M\ot\pno(M\ot M)\ar[dl]\ar[d] &
\text{4 factors,}\atop\text{3 arrows}\ar[dl]\ar[d]\\
\pno M& \pno M\ot\pno M & \pno M\ot\pno M\ot \pno M &
         \pno M\pno M\ot\pno M \ot\pno M}$$
For $i\leq 0$, the degree $i$ component $\Omega\B A^i$ 
of $\Omega\B A$ is the sum
over all terms of the form $\pno(M\ots M)\ots \pno(M\ots M)$,
containing $n=1,2,3,\ldots$ copies of $M$ and $n+i$ arrows.
We call the number $n$ the \textbf{polynomial degree}.
The polynomial degree is a second grading of the algebra $\Omega\B A$.
A third one is given by the number of arrows arising in a summand
The three gradings are related by
$$\text{degree}=\text{number of arrows}-\text{polynomial degree}.$$
The differential $d$ of $\Omega\B A$ has two components.
One component $d_{sw}$ of polynomial degree $-1$ and one component
$d_s$ increasing the number of arrows by one.
Since $\Omega\B A$ is the free tensor algebra over $\pno(TM)$,
the differential is uniquely defined by its values on
terms of the form $\pno(M\ots M)$ (upper diagonal).
We have 
$$d_s(\pno(m_1\ots m_n))=\sum_{j=1}^{n-1}(-1)^j\pno(m_1\ots m_j)\ot\pno
(m_{j+1}\ots m_n)$$
and $d_{sw}$ is just the shifted bar complex, i.e.
$$d_{sw}(\pno(m_1\ots m_n))=\sum_{j=1}^{n-1}(-1)^{j-1}
\pno(m_1\ots m_{j-1}\ot m_j\cdot m_{j+1}\ots m_n),$$
where $m_j\cdot m_{j+1}$ stands for $\pnu(\pno m_j\cdot\pno m_{j+1})$.
The natural algebra map $p:\Omega\B A\nach A$ given 
by the adjoint functor property
of the bar and cobar construction is the map from the degree
zero componet to $A$ given by
$$p:\pno m_1\ots \pno m_n\mapsto \pno m_1\ppp \pno m_n.$$
 
The complex $\Omega \B A$ is filtered by the polynomial degree
and the associated differential graded ring
is the sum over the rows in the above picture.
All rows except the first one are acyclic split.
The splitting is given by the formular (see 
Chapter 3.4 of \cite{Smir} or Page 36 of \cite{Lefe})
\begin{equation*}
\pno(m_{11}\ots m_{1i_1})\ots \pno(m_{k1}\ots m_{ki_k})\mapsto
\end{equation*}
\begin{equation*}
\left\{ \begin{array}{r@{\quad\text{ for }\quad}l}
0 & i_1\neq 1\\
\pno(m_{11}\ot m_{21}\ots m_{2i_2})\ots\pno(m_{k1}\ots m_{ki_k}) 
& i_1=1
\end{array}\right. 
\end{equation*}
Again by convergence of the spectral sequence obtained from
the polynomial filtration,
the map $p$ is a quasi-isomorphism. By the same reason,
the natural inclusion 
$A=\pno M\nach\Omega \B A$ is also a quasi-isomorphism. Of course it
is not an algebra map. The proof of the following theorem is
encrypted in \cite{KontL}:
\begin{satz}\label{KKK}
There is a natural quasi-isomorphism
$$L_{\rel}=\Der(\Omega \B A)\nach \Coder(\B A)=L_{\ass}$$
respecting the DGL-structures. In consequence,
we get an isomorphism $\Def_{rel}\nach\Def_{\ass}$
of deformation functors.
\end{satz}
\begin{proof}
The quasi-isomorphism $p:\Omega \B A\nach A$ induces a quasi-isomorphism
$$\Der(\Omega\B A)=\Hom(\pno\B A,\Omega\B A)\nach\Hom(\pno\B A,A)\isom
\Hom(\B A,\pnu A)=\Coder(\B A).$$
Check that it respects the DGL structures.
\end{proof}

In the sequel, we will form the set $\MC_{\rel}(B)$ with respect to
the canonical free DGA resolution $\Omega\B A$ of $A$.
By Theorem~\ref{isodef}, the quasi-isomorphism of Theorem~\ref{KKK} 
induces, for each $B\in\Art$, a bijection
$E:\Def_{\rel}(B)\nach\Def_{\ass}(B)$. The knowledge of the
structure of $\Omega\B A$ allows us to describe easily
the inverse $E^{-1}$: Let $\beta\in\MC_{\ass}(B)$ be
a Maurer-Cartan element. The differentials of
$\Omega\B(A_B,\alpha)$ and $\Omega\B(A_B,\alpha+\beta)$
on the graded algebra $\Omega\B A_B$ 
only differ by a component $\delta(\beta)$ in south-west direction.
We have $\delta(\beta)\in\MC_{\rel}(B)$ and the quasi-isomorphism
of Theorem~\ref{KKK} maps the
class of $\delta(\beta)$ to $\beta$.
Thus, $E$ maps the class $[\delta(\beta)]$ of $\delta(\beta)$  
in $\Def_{\rel}(B)$ to the class $[\beta]$ of
$\beta$ in $\Def_{\ass}(B)$. 

Now we are able to prove the main statement:

\begin{satz}\label{drei}
For any local commutative artinian algebra $B$,
the diagram 
$$\xymatrix{
\Def_{\rel}(B)\ar[r]^E\ar[dr]_G & \Def_{\ass}(B)\ar[d]^F\\
 & \Def_{\flat}(B)
}$$
commutes and all the arrows are bijections.
\end{satz}
\begin{proof}
For $[\beta]\in\Def_{\ass}(B)$, we have
$G(E^{-1}([\beta]))=G([\delta(\beta)])=H^0(\Omega\B A_B,d+\delta(\beta))\isom
(A_B,\alpha+\beta)\isom F([\beta])$.
This proves the commuativity.
Since $E$ is an bijective, $G$ is surjective and $F$ is injective,
we see that $F$ and $G$ are also isomorphisms. 
\end{proof}

\begin{kor}\label{vier}
The three deformation functors $\Def_{\ass}$, $\Def_{\rel}$ and
$\Def_{\flat}$ are isomorphic.
\end{kor}

By Theorem~\ref{drei}, the functors $E$, $F$ and $G$ are
in particular smooth. Thus, by an observation of Schlessinger
\cite{Schl}, Corollary~\ref{vier} is still true, if we
replace the three deformation functors by their extensions
to certain categories of projective limits of Artin algebras,
for example the category of complete noetherian local
algebras with residue field $k$.

\end{document}